\documentclass[12pt]{amsart}
\usepackage[english]{babel}
\usepackage{amssymb,amsmath,amscd,amsthm}

\usepackage{cite}
\def\?{}
\newtheoremstyle{neosn}{0.5\topsep}{0.5\topsep}{\rm}{}{\sc}{.}{ }{\thmname{#1}\thmnumber{ #2}\thmnote{ {\mdseries#3}}}
\theoremstyle{neosn}

\newcommand{\Soc}{\,\mathrm{Soc}\,}

\begin{document}
\begin{center}
\Large{\textbf{Rings Which Are Essential over Their Centers}}
\end{center}

\hfill {\sf V.T. Markov}

\hfill Lomonosov Moscow State University

\hfill e-mail: vtmarkov@yandex.ru

\hfill {\sf A.A. Tuganbaev}

\hfill National Research University "MPEI"

\hfill Lomonosov Moscow State University

\hfill e-mail: tuganbaev@gmail.com

{\bf Abstract.} A centrally essential ring is a ring which is an essential extension of its center (we consider the ring as a module over its center). We give several examples of noncommutative centrally essential rings and describe some properties of centrally essential rings.

V.T.Markov is supported by the Russian Foundation for Basic Research, project 17-01-00895-A. A.A.Tuganbaev is supported by Russian Scientific Foundation, project 16-11-10013.

{\bf Key words:}
centrally essential ring, center of the ring, semiprime ring, semiperfect ring

\begin{center}
\textbf{1. Introduction}
\end{center}
All considered rings are associative unital rings.

A ring $R$ with center $C$ is said to be \emph{centrally essential} if the module $R_C$ is an essential extension of the module $C_C$.

It is clear that any commutative ring is centrally essential. The following example shows that the converse is not always true even for finite rings.

\textbf{Example.} Let $F$ be the field consisting of three elements, $V$ be a vector $F$-space with basis $e_1,e_2,e_3$, and let $\Lambda(V)$ be the exterior algebra of the space $V$, i.e. $\Lambda(V)$ is the algebra with operation $\wedge$ defined by generators $e_1,e_2,e_3$ and defining relations
$$
e_i\wedge e_j+e_j\wedge e_i=0\mbox{ for all } i,j=1,2,3.
$$
It is clear that $e_1^2=e_2^2=e_3^2=0$ and any product of
generators is equal to the $\pm$product of generators with
ascending subscripts.

Thus,  $\Lambda(V)$ is a finite $F$-algebra of dimension 8 with basis
$$
\{1,e_1,e_2,e_3,e_1\wedge e_2,e_1\wedge e_3,e_2\wedge e_3,e_1\wedge e_2\wedge e_3\}
$$
and $|\Lambda(V)|=3^8$.

It is clear that $e_k\wedge e_i\wedge e_j=-e_i\wedge e_k\wedge e_j=e_i\wedge
e_j\wedge e_k$. Therefore, if
$$x=\alpha_0\cdot
1+\alpha_1^1e_1+\alpha_1^2e_2+\alpha_1^3e_3
+\alpha_2^1e_1\wedge
e_2+\alpha_2^2e_1\wedge e_3+
\alpha_2^3e_2\wedge e_3+\alpha_3e_1\wedge e_2\wedge
e_3,$$ then
{\renewcommand{\arraycolsep}{0pt}
$$
\begin{array}{lll}[e_1,x]&=&2\alpha_1^2e_1\wedge
e_2+2\alpha_1^3e_1\wedge e_3,\\
{}[e_2,x]&=-&2\alpha_1^1e_1\wedge
e_2+2\alpha_1^3e_2\wedge e_3,\\
{}[e_3,x]&=-&2\alpha_1^2e_1\wedge
e_3-2\alpha_1^2e_2\wedge
e_3.\end{array}$$ } Thus, $x\in C(\Lambda(V))$ if and only if  $\alpha_1^1=\alpha_1^2=\alpha_1^3=0.$ In other words, the center of the algebra $\Lambda(V)$ is of dimension 5. On the other hand,
if, e.g., $\alpha_1^1\neq 0$, then
$$
x\wedge(e_2\wedge e_3)=\alpha_0e_2\wedge e_3+\alpha_1^1e_1\wedge e_2\wedge e_3\in C(\Lambda(V))\setminus \{0\}.$$ In addition,
$e_2\wedge e_3\in C(\Lambda(V))$. Therefore, $\Lambda(V)$ is a
finite centrally essential noncommutative ring.

The main results of the paper are Theorem 1, 2 and 3.

\textbf{Theorem 1.} \textit{Let $R$ be a ring.}

\textbf{1.} \textit{If $R=\oplus_{n\in\mathbb{N}_0}R_n$ is a graded generalized anticommutative homogeneously faithful ring $($see the definitions below$)$ without additively 2-torsion elements,
then the ring $R$ is centrally essential if and only if either $R=R_0$ or
there exists an odd positive integer $n$ such that $R_n\neq 0$ and $R_{n+1}=0$.}

\textbf{2.} \textit{The exterior algebra $\Lambda(V)$ of a finite-dimensional vector space $V$ over a field $F$ of characteristic $0$ or $p\neq 2$ is a centrally essential ring if and only if
$\dim V$ is an odd positive integer. In particular, if $F$ is a finite field of odd characteristic and $\dim V$ is an odd positive integer exceeding $1$, then $\Lambda(V)$ is a centrally essential noncommutative finite ring.}

\textbf{3.} \textit{There exist noncommutative centrally essential finite rings.}

\textbf{Theorem 2.} \textit{Let $R$ be a centrally essential ring with center $C$.}

\textbf{1.} \textit{If $C$ is a semiprime ring, then the ring $R$ is commutative.}

\textbf{2.} \textit{If the ring $R$ is semiprime, then $R$ is commutative.}

\textbf{3.} \textit{In the ring $R$, all idempotents are central.}

\textbf{4.} \textit{If the ring $R$ is semiperfect, then $R/J(R)$ is a finite direct product of fields; in particular, $R/J(R)$ is a
commutative ring. In addition, $R$ is a finite direct product of centrally essential local rings and $\Soc(R_C)\subseteq C$.}

\textbf{Theorem 3.} \textit{If $R$ is a right or left perfect ring with center $C$, then the following conditions are equivalent.}
\begin{enumerate}
\item[\textbf{1)}]
$R$ is a \textit{centrally essential ring.}
\item[\textbf{2)}]
$\Soc(R_C)\subseteq C$ \textit{and all idempotents of the ring $R$ are central.}
\end{enumerate}

The proof of Theorem 1, Theorem 2, and Theorem 3 is decomposed into a series of assertions, some of which are of independent interest.
We give some necessary notions. All remaining necessary notions of ring theory can be found in the books\cite{Bourbaki,Lambek,Tuganbaev,Herstein,Lam}.

We denote by $C(R)$ and $J(R)$ the center and the Jacobson radical of the ring $R$, respectively; we also set $[a,b]=ab-ba$ for any two elements $a,b$ of the ring $R$.

A ring $R$ is said to be \emph{local}, if $R/J(R)$ is a division ring.

A ring $R$ is said to be \emph{semiperfect} if $R/J(R)$ is an
Artinian ring and every complete system of orthogonal idempotents of the ring $R/J(R)$ can be lifted to a complete system of orthogonal idempotents of $R$.

A ring $R$ is said to be \textit{left perfect} if $R$ is semiperfect and its Jacobson radical $J(R)$ is left $T$-\textit{nilpotent}, i.e., for any sequence  $x_1,x_2,\ldots$ of elements of $J(R)$,
there exists a subscript $n$ such that $x_1x_2\ldots x_n=0$. Right perfect rings are similarly defined.

For a module $M$, the \emph{socle} $\Soc M$ of $M$ is the sum of all simple submodules in $M$; if $M$ does not have a simple submodule, then $\Soc M=0$ by definition.

Let $(S,+)$ be a semigroup. A ring $R$ is said to be $S$-graded if
$R$ is the direct sum of additive subgroups $R_s,\ s\in S$ and
$R_sR_t\subseteq R_{s+t}$ for any elements $s,t\in S$. For any $s\in S$, the elements of the subgroup $R_s$ are called \textit{homogeneous} elements of degree $s$.

If $S=\mathbb{N}_0=\mathbb{N}\cup\{0\}$, then $S$-graded rings are called \textit{graded} rings. It is easy to verify that the identity element of a graded ring is contained in the subgroup $R_0$. On an
arbitrary graded ring $R=\oplus_{n\in\mathbb{N}_0}R_n$, we can define a $\mathbb{Z}_2$-grading:\\
$R=R_{(0)}\oplus R_{(1)}$, where $R_{(i)}=\bigoplus_{k\in
\mathbb{N}_0}R_{2k+i}$, $i\in \{0,1\}.$

We say that the graded ring $R=\oplus_{n\in\mathbb{N}_0}R_n$ is \emph{generalized anticommutative} if for any two integers $m,n\in\mathbb{N}_0$ and every elements $x\in R_m$ and $y\in R_n$, the relation $yx=(-1)^{mn}xy$ holds.

If the graded ring $R=\oplus_{n\in\mathbb{N}_0}R_n$ satisfies the condition
$$
\forall m,n\in \mathbb{N}_0,\quad \ R_{m+n}\neq 0\ \Rightarrow\ R_m\neq 0\ \&\ \forall x\in{R_m}\setminus\{0\}, xR_n\neq 0,\eqno{(*)}
$$ then we say that $R$ is a \emph{homogeneously faithful} ring.

Let $F$ be a field of characteristic $0$ or $p\neq 2$, $V=F^n$ be a vector space over $F$ of dimension $n>0$, and let $\Lambda(V)$ be the exterior algebra of the space $V$ \cite[\S III.5]{Bourbaki}, which can be defined as a unital $F$-algebra with respect to multiplication operation $\wedge$ with generators $e_1,\ldots,e_n$ and defining relations $e_i\wedge e_j+e_j\wedge e_i=0$ for all $i,j\in\{1,\ldots,n\}$.

The algebra $\Lambda(V)$ has a natural grading:\\
$\Lambda(V)=\bigoplus_{p\in\mathbb{N}_0}\Lambda^p(V)$,
where $\Lambda^p(V)$ for $1\leq p\leq n$, is a vector space
with basis
$$
\{e_{i_1}\wedge \ldots \wedge e_{i_p}:1\leq i_1<\ldots<i_p\leq
n\},
$$
$\Lambda^0(V)=F$ and $\Lambda^p(V)=0$ for $p>n$.

It is well known that exterior algebras are generalized
anticommutative.

\begin{center}
\textbf{2. Graded Centrally Essential Rings}
\end{center}

\textbf{Proposition 2.1.}\label{center} \textit{In an arbitrary graded ring $R=\oplus_{n\in\mathbb{N}_0}R_n$, the following relation $C(R)=\oplus_{n\in\mathbb{N}_0}(R_n\cap C(R))$ holds}.

\textbf{Proof.} The inclusion $\oplus_{n\in\mathbb{N}_0}(R_n\cap C(R))\subseteq C(R)$ is obvious. Let $x=x_0+x_1+\ldots x_n\in C(R)$, where $x_i\in R_i$, $i=0,1,\ldots,n$. If $y\in R_m$, for some $m\in\mathbb{N}_0$, then $0=[x,y]=[x_0,y]+\ldots+[x_n,y]$ and summands of the last sum are contained in different direct summands $R_{m},R_{m+1},\ldots,R_{m+n}$. Therefore, $[x_i,y]=0$ for each homogeneous element $y$ and any $i=0,1,\ldots,n$. Then $x_i\in C(R)$, since any element of the ring is the sum of homogeneous elements.~\hfill$\square$

\textbf{Remark.} The proof of Proposition 2.1 remains true for each $S$-graded ring provided $S$ is a commutative cancellative semigroup.

\textbf{Proposition 2.2.}\label{exter}
\textit{The graded algebra $R=\Lambda(V)$ is a homogeneously faithful ring.}

\textbf{Proof.}
Let $p,q\in\{0,\ldots,n\}$ and $p+q\leq n$. If $pq=0$, then the condition $(*)$ holds. Now let $0<p<n$ and $0\neq x\in R_p$. We take
a basis element $e_{i_1}\wedge \ldots \wedge e_{i_p}$ which has the non-zero coefficient in the representation of $x$. Since
$p+q\leq n$, there exist subscripts $j_1,\ldots,j_q\in\{1,\ldots,n\}$ such that $1\leq j_1< \ldots <j_q\leq
n$ and $\{i_1,\ldots,i_p\}\cap\{j_1,\ldots,j_q\}=\varnothing$.
We set $y=e_{j_1}\wedge\ldots\wedge e_{j_q}$ and note that the basis element $\pm e_{i_1}\wedge \ldots \wedge e_{i_p}\wedge
e_{j_1}\wedge\ldots\wedge e_{j_q}$ of the space
$\Lambda^{p+q}(V)$ has the non-zero coefficient in the representation of the element $xy$, since the products of remaining basis elements of the space $\Lambda^{p}(V)$ by the element $y$ either are equal to $0$ or are equal to $\pm$ other basis elements of the space $\Lambda^{p+q}(V)$.~\hfill$\square$

\textbf{Proposition 2.3.}\label{centerdescr} \textit{Let $R=\oplus_{n\in\mathbb{N}_0}R_n$ be a graded generalized anticommutative homogeneously faithful ring which does not have additively 2-torsion elements. If there exists an odd positive integer $n$ such that $R_n\neq 0$ and $R_{n+1}=0$, then $C(R)=R_{(0)}+R_n$. Otherwise, $C(R)=R_{(0)}$.}

{\bf Proof.} It follows from the relation of generalized anticommutativity that $R_{(0)}\subseteq C(R)$. The following property follows from $(*)$:\\
if such an integer $n$ exists, then $R_m=0$ for $m>n$ and $R_m\neq 0$ for $0\leq m \leq n$, in addition, if $x\in R_n$ and
$y=y_0+z\in R$, where $y_0\in R_0$ and $z\in\oplus_{m>0}R_m$, then
$[x,y]=[x,y_0]=0$, i.e., $R_n\subseteq C(R)$. Conversely, let $x\in
C(R)$. By Proposition 2.1, we can assume that $x$ is a homogeneous element of odd degree $i$. Let $x \neq 0$ and $R_{i+1}\neq 0$. Then it follows from $(*)$ that there exists an element $y\in R_1$ with $xy\neq 0$. We obtain $0=[x,y]=2xy$; this is a contradiction. Thus, either $x=0$ or $x\neq 0$, but $R_{i+1}=0$, i.e., $i=n$.~\hfill$\square$

\textbf{Proposition 2.4.}\label{esscenter}
\textit{Let $R=\oplus_{n\in\mathbb{N}_0}R_n$ be a graded generalized
anticommutative homogeneously faithful ring without additively
2-torsion elements. The ring $R$ is centrally essential if and
only if either $R=R_0$ or there exists an odd positive integer $n$
such that
$R_n\neq 0$ and $R_{n+1}=0$.}

{\bf Proof.} Let $R$ be a centrally essential ring, $C=C(R)$ and $R\neq R_0$. By $(*)$, we have $R_1\neq 0$. We take an element $x\in R_1\setminus\{0\}$ and assume that such an integer $n$ does not exist. By Proposition 2.3, we have $C=R_{(0)}$ and $xC\subseteq R_{(1)}$, whence $xC\cap C\subseteq R_{(1)}\cap R_{(0)}=0$; it is a contradiction.

Conversely, if $R=R_0$, then $C=R$, since the ring $R_0$ is commutative. We assume that there exists an odd positive integer $n$ such that $R_n\neq 0$ and $R_{n+1}=0$. Let $0\neq x\in R\setminus C$. We have $x=x_0+\ldots+x_n$, where $x_i\in R_i$, and we take the least odd positive integer $m$ with $x_m\neq 0$. It is clear that $1\leq m\leq n$. We set $k=n-m$ and take an element $y\in R_k$ such that $x_my\neq 0$. It is clear that $y\in C$. In addition, $xy$ is the sum of homogeneous elements of even degree and the element $x_my$ of odd degree $n$. Therefore, $xy\in C$, by Proposition 2.3, and $xy\neq 0$.~\hfill$\square$

The following proposition follows from Proposition 2.2 and Proposition 2.4.

\textbf{Proposition 2.5.}
\textit{Let $V$ be a finite-dimensional vector space over a field $F$ of characteristic $0$ or $p\neq 2$. The exterior algebra $\Lambda(V)$ of $V$ is a centrally essential ring if and only if $V$ is of odd dimension. In particular, if $F$ is a finite field of odd characteristic and $\dim V$ is an odd positive integer exceeding $1$, then $\Lambda(V)$ is a centrally essential noncommutative ring.}

\begin{center}
\textbf{3. Quasi-Identity and Idempotents\\ of a Centrally
Essential Ring}
\end{center}

\textbf{Proposition 3.1.}\label{quid}
Any centrally essential ring $R$ satisfies the following relations.
\begin{multline}\label{eqquid}
\forall n\in\mathbb{N}, x_1,\ldots,
x_n,y_1,\ldots,y_n,r\in R,\\
\left\{\begin{array}{rcl}x_1y_1+\ldots+x_ny_n&=&1\\
x_1ry_1+\ldots+x_nry_n&=&0\end{array}\right.
\Rightarrow r=0.
\end{multline}

{\bf Proof.} We assume that the above relations hold, but $r\ne 0$. Then there exist two elements $c,d\in C(R)$ such that $cr=d\neq 0$. Consequently,
$$
d=d(x_1y_1+\ldots+x_ny_n)=x_1dy_1+\ldots+x_ndy_n=c(x_1ry_1+\ldots+x_nry_n)=0.
$$
This is a contradiction.~\hfill$\square$

\textbf{Proposition 3.2.}\label{idemp}
\textit{If $R$ is a centrally essential ring, then all idempotents  of $R$ are central.}

{\bf Proof.} Let $e\in R$ and $e^2=e$. Then $e^2+(1-e)^2=1$. For each $x\in R$, we have $e(ex-exe)e+(1-e)(ex-exe)(1-e)=0$ and
$e(xe-exe)e+(1-e)(xe-exe)(1-e)=0$. By Proposition 3.1, we have
$xe-exe=ex-exe=0$, whence $xe=exe=ex$.~\hfill$\square$

\textbf{Proposition 3.3.}\label{loc}
\textit{If $R$ is a centrally essential local ring, then the ring $R/J(R)$ is commutative.}

{\bf Proof.} Let $x,y\in R$ and $xy-yx\not \in J(R)$. If $c\in C(R)$ and $cx=d\in C(R)\setminus \{0\}$, then $c(xy-yx)=dy-yd=0$, whence $c=0$. This is a contradiction.~\hfill$\square$

\textbf{Proposition 3.4.}\label{semiperf} \textit{Let $R$ be a centrally essential semiperfect ring. Then $R/J(R)$ is a finite direct product of fields; in particular, $R/J(R)$ is a commutative ring. In addition, $R$ is a finite direct product of centrally essential local rings.}

{\bf Proof.} By the definition of a semiperfect ring, $R/J(R)$ is the~direct sum of simple Artinian rings each of them is the matrix ring over a division ring. Let $\bar{e}_1,\ldots,\bar{e}_n$ be a complete system of indecomposable orthogonal idempotents of
$\bar R=R/J(R)$. Then there exists a complete system of indecomposable orthogonal idempotents $e_1,\ldots,e_n$ in $R$ such that $e_i+J(R)=\bar{e}_i$, $i=1,\ldots,n$. By Proposition 3.2, all the idempotents $e_1,\ldots,e_n$ are central. Therefore, $R=\oplus_{i=1}^n R_ie_i$ is a decomposition of the ring $R$ into the direct sum of local centrally essential rings. Consequently, all the rings $R_i/J(R_i)$ are commutative by Proposition 3.3. It is easy to verify that all the rings $R_i=Re_i$ are centrally essential; therefore, the division ring $R_i/J(R_i)$ is commutative. Then the ring $R/J(R)=\oplus_{i=1}^n R_i/J(R_i)$ is also commutative.~\hfill$\square$

\textbf{Proposition 3.5.}\label{soccenter}
\textit{If $R$ is a centrally essential semiperfect ring with center $C$, then $\Soc(R_C)\subseteq C$.}

\textbf{Proof.} It follows from the proof of Proposition 3.4 that we can assume, without loss of generality, that $R$ is a local ring.

We remark that $C=C(R)$ is a local ring and $J(C)=C\cap J(R)$.

Now let $x\in\Soc(R_C)\setminus\{0\}$. There exist two elements $c,d\in C$ with $cx=d\neq 0$. It is clear that $c\not\in J(R)$, since
$J(C)\Soc(R_C)=0$. Consequently, $c$ is an invertible element and $x=c^{-1}d\in C$.~\hfill$\square$

Obviously, it follows from the above that $\Soc({}_RR)=\Soc(R_R)$ if $R$ is a centrally essential semiperfect ring.

\begin{center}
\textbf{4. The Completion of the Proof of Theorem 1, Theorem 2, and Theorem 3.}
\end{center}

\textbf{Proposition 4.1.}\label{semiprime}
\textit{If $R$ is a centrally essential ring and $C=C(R)$ is a
semiprime ring, then the ring $R$ is commutative.}

\textbf{Proof.} For each $r\in R$, we verify that the ideal $$r^{-1}C=\{c\in C:\;rc\in C\}$$ is dense in $C$. Indeed, let $d\in C$ and
$dr^{-1}C=0$. If $dr=0$, then $d\in r^{-1}C$ and $d^2=0$, whence
$d=0$. Otherwise, since $R$ is a centrally essential ring, there exists an element $z\in C$ with $zdr\in C\setminus\{0\}$. Then $zd\in r^{-1}C$ and $(zd)^2=0$, whence $zd=0$; this contradicts to the choice of $z$. The property `$r(r^{-1}C)\neq 0$ for each $r\in R\setminus\{0\}$' is equivalent to the property that the ring $R$ is centrally essential. Therefore, $R$ is a right ring of quotients of the ring $C$ in the sense of \cite[\S 4.3]{Lambek}; therefore, $R$ can be embedded in the complete ring of quotients of the ring $C$ which is commutative.~\hfill$\square$

\textbf{4.2. The completion of the proof of Theorem 1.}

Theorem 1 follows from Proposition 2.4 and Proposition 2.5.

\textbf{4.3. The completion of the proof of Theorem 2.}

\textbf{1.} The assertion is proved in Proposition 4.1.

\textbf{2.} Let $R$ be a centrally essential semiprime ring with center $C$. By Proposition 4.1, it is sufficient to prove that $C$ is a semiprime ring. Let $c\in C$ and $c^2=0$. Since $C=C(R)$, we have $(RcR)^2=0$. Since the ring $R$ is semiprime, $c=0$ and the
ring $C$ is semiprime.

\textbf{3.} The assertion is proved in Proposition 3.2.

\textbf{4.} Let $R$ be a centrally essential semiperfect ring. By Proposition 3.4, $R/J(R)$ is a finite direct product of fields and, in particular, a commutative ring; in addition, $R$ is a finite direct product of centrally essential local rings. By Proposition 3.5, $\Soc(R_C)\subseteq C$.

\textbf{4.4. The completion of the proof of Theorem 3.}

1)\,$\Rightarrow$\,2). The assertion follows from Proposition 3.2 and Proposition 3.5.

2)\,$\Rightarrow$\,1). Again, since all idempotents are central, we can assume that $R$ is a local ring. Then $J(C)=C\cap J(R)$ and $C/J(C)$ is a field.

Let $x\in R\setminus \{0\}$. If $J(C)x=0$, then $x\in\Soc(R_C)$; therefore, $x\in C$. Otherwise, there exists an element $c_1\in J(C)$ such that $c_1x\neq 0$. If $J(C)c_1x=0$, then $c_1x\in\Soc(R_C)$ and $c_1x\in C$; otherwise, we take an element $c_2\in J(C)$ with $c_2c_1x\neq 0$ and so on. Since the radical $J(R)$ of a right or left perfect ring $R$ is right or left $T$-nilpotent and the elements $c_i$ are central, this process will stopped at some finite step.~\hfill$\square$

\textbf{Remark.} In Theorem 3, we cannot omit the condition that the ring $R$ is right or left perfect, since every noncommutative local ring
without zero-divisors, e.g., the formal power series ring in one
variable over the division ring of Hamiltonian quaternions,
satisfies all remaining conditions of this theorem, but this ring
is not centrally essential.

\textbf{4.5. Open questions.} Let $R$ be a centrally essential ring with center $C$.

\textbf{1.} Is it true that the factor ring of the ring $R$ with respect to its Jacobson radical or its prime radical is commutative?

\textbf{2.} If the ring $R$ is semiperfect, is it true that $\Soc(R_C)=\Soc(R_R)$?

\textbf{3.} If the ring $R$ is semiperfect, is it true that $R=C+J(R)$?

\end{document}